\documentstyle[amstex,amssymb, 12pt]{article}

\renewcommand{\phi}{\varphi}

\newtheorem{theorem}{Theorem}[section]

\newtheorem{corollary}[theorem]{Corollary}

\newtheorem{definition}[theorem]{Definition}

\begin{document}
\thispagestyle{empty}

\begin{center}
{\Large\bf Characterizing Type I $C^*$-algebras via Entropy}

\end{center}

\begin{center}{\bf Nathanial P. Brown \footnote{ 
      Currently an  MSRI Postdoctoral Fellow. 

      1991 AMS Classification; 46L55.}}\\
   UC-Berkeley\\
   Berkeley, California 94720 \\
   nbrown{\char'100}math.berkeley.edu 

\end{center}

\begin{abstract}
Let $A$ be a separable unital $C^*$-algebra.  It is shown that $A$ is
type I if and only if the CNT-entropy of every inner automorphism of
$A$ is zero.
\end{abstract}

\parskip2mm

\section{Introduction}

The theory of noncommutative entropy began in the context of
noncommutative dynamical systems.  However, in this short note we
observe that entropy is also closely connected to the structure theory
of $C^*$-algebras.  Taking this point of view it will be useful to
regard entropy as an invariant of an algebra, as opposed to an
invariant of a dynamical system.  Consider the following analogue of
the topological entropy invariants discussed in [BDS, 6.4.1].

\begin{definition}
Let $A$ be a unital $C^*$-algebra.  Let
$CNT_{Inn} (A)$ denote the set of real numbers $t$ such that there
exists a unitary $u \in A$ and a state $\phi \in S(A)$ such that
$\phi \circ {\rm Ad}u = \phi$ and $h_{\phi}({\rm Ad}u) = t$,
where $h_{\phi}({\rm Ad}u)$ denotes the CNT-entropy [CNT].
\end{definition}

In [NS] it is shown that the CNT-entropy in a type I $W^*$-algebra
can always be computed by looking at the restriction to the center.
Since inner automorphisms act trivially on the center this implies that
{\em if $A$ is a type I $C^*$-algebra then $CNT_{Inn} (A) = \{0\}$} (see
[NS, Cor.\ 7]).  Hence this entropy invariant provides a natural
dynamical obstruction to an algebra being type I.  Our main result
shows that this is the only obstruction, i.e. we observe the converse
of Neshveyev and St$\o$rmer's result.

\begin{theorem}
Let $A$ be a separable unital $C^*$-algebra.  Then $A$ is type I if
and only if $CNT_{Inn} (A) = \{0\}$.
\end{theorem}

{\noindent\bf Proof.}  We refer to [CNT] or [St] for all relevant
definitions and the notation which follows. 

In light of [NS, Cor.\ 7], it suffices to show that if $A$ is unital,
separable and {\em not} type I then there exists a unitary $u \in A$
and a state $\phi \in S(A)$ such that $\phi \circ {\rm Ad}u = \phi$
and $h_{\phi}({\rm Ad}u) > 0$.  This will follow from Glimm's theorem
and a few other well known results.

So assume $A$ is unital, separable and {\em not} type I.  Let ${\cal
U} = \otimes_{n \in {\Bbb N}} M_n ({\Bbb C})$, $M_{2^{\infty}} =
\otimes_{\Bbb Z} M_2 ({\Bbb C})$ be the CAR algebra and $\gamma \in
Aut(M_{2^{\infty}})$ be the noncommutative Bernouli shift on
$M_{2^{\infty}}$ arising from the map ${\Bbb Z} \to {\Bbb Z}$, $i
\mapsto i + 1$.  By [Vo] we may regard the crossed product
$M_{2^{\infty}} \rtimes_{\gamma} {\Bbb Z}$ as a unital subalgebra of
${\cal U}$.  By Glimm's theorem [Pe, Thm.\ 6.7.3] we can find a unital
subalgebra $\tilde{B} \subset A$ and a surjective $*-$homomorphism
$\pi : \tilde{B} \to {\cal U}$.  Let $B = \pi^{-1}(M_{2^{\infty}}
\rtimes_{\gamma} {\Bbb Z}) \subset \tilde{B} \subset A$.  Since the
unitary group of ${\cal U}$ is connected, any unitary in ${\cal U}$
lifts to a unitary in $\tilde{B}$.  So, letting $v \in M_{2^{\infty}}
\rtimes_{\gamma} {\Bbb Z}$ be the implementing unitary, we can find a
unitary $u \in B$ such that $\pi(u) = v$. 

It should now be clear what to do: use the trace on $M_{2^{\infty}}
\rtimes_{\gamma} {\Bbb Z}$ to construct an ${\rm Ad}u$-invariant state
on $A$ with positive entropy.  This requires a bit of care, but is
essentially straightforward.

Let $\tau$ be the unique tracial state on $M_{2^{\infty}}
\rtimes_{\gamma} {\Bbb Z}$ and $(\pi_{\tau}, L^2( M_{2^{\infty}}
\rtimes_{\gamma} {\Bbb Z}, \tau), \eta_{\tau})$ be the associated GNS
representation (here, $\eta_{\tau}$ denotes the canonical cyclic
vector in $L^2( M_{2^{\infty}} \rtimes_{\gamma} {\Bbb Z}, \tau)$).
Since $M_{2^{\infty}} \rtimes_{\gamma} {\Bbb Z}$ is nuclear, we can
apply Arveson's extension theorem followed by a conditional
expectation to construct a unital completely positive map $\Phi : A
\to \pi_{\tau} (M_{2^{\infty}} \rtimes_{\gamma} {\Bbb
Z})^{\prime\prime} \subset B(L^2( M_{2^{\infty}} \rtimes_{\gamma}
{\Bbb Z}, \tau))$ which is an extension of the map $\pi_{\tau} \circ
\pi : B \to \pi_{\tau} (M_{2^{\infty}} \rtimes_{\gamma} {\Bbb Z})$.
Define a state $\phi$ on $A$ by $$\phi(a) = <\Phi(a) \eta_{\tau},
\eta_{\tau}>.$$

Since $\Phi$ maps $u$ to a unitary operator, $u$ lies in the
multiplicative domain of $\Phi$ (see [Pa, Exercise 4.2]).  Because we
have arranged that $\Phi$ takes values in $\pi_{\tau} (M_{2^{\infty}}
\rtimes_{\gamma} {\Bbb Z})^{\prime\prime}$, this will imply that $\phi
\circ {\rm Ad}u = \phi$.  To see this we first observe that
$\pi_{\tau} (v) T \pi_{\tau} (v)^* = U_{\gamma} T U_{\gamma}^*$, for
all $T \in \pi_{\tau} (M_{2^{\infty}} \rtimes_{\gamma} {\Bbb
Z})^{\prime\prime}$ where $U_{\gamma} : L^2( M_{2^{\infty}}
\rtimes_{\gamma} {\Bbb Z}, \tau) \to L^2( M_{2^{\infty}}
\rtimes_{\gamma} {\Bbb Z}, \tau)$ is the unitary operator defined by
$U_{\gamma} (x) = vxv^*$.  Hence we have \begin{multline*}
\begin{aligned} \phi(uau^*)
    &=  <\Phi(uau*) \eta_{\tau}, \eta_{\tau}> \\[2mm]
    &=  <\pi_{\tau}(v) \Phi(a) \pi_{\tau}(v)^* \eta_{\tau}, \eta_{\tau}>\\[2mm]
    &=  <U_{\gamma} \Phi(a) U_{\gamma}^* \eta_{\tau}, \eta_{\tau}>\\[2mm] 
    &=  \phi(a), 
\end{aligned}
\end{multline*}
for all $a \in A$. 

To finish the proof we must observe that $h_{\phi} ({\rm Ad}u) \geq
h_{\tau} ({\rm Ad}v) = \log 2$.  This almost follows immediately from
the definitions.  The only observation which needs to be made is that
the CNT-entropy of the systems $(\pi_{\tau} (M_{2^{\infty}}
\rtimes_{\gamma} {\Bbb Z}), {\rm Ad}\pi_{\tau}(v), < \cdot
\eta_{\tau}, \eta_{\tau}>)$ and $(\pi_{\tau} (M_{2^{\infty}}
\rtimes_{\gamma} {\Bbb Z})^{\prime\prime}, {\rm Ad}\pi_{\tau}(v), <
\cdot \eta_{\tau}, \eta_{\tau}>)$ naturally agree.  Hence when
computing the CNT-entropy of the dynamical system $(\pi_{\tau}
(M_{2^{\infty}} \rtimes_{\gamma} {\Bbb Z})^{\prime\prime}, {\rm
Ad}\pi_{\tau}(v), < \cdot \eta_{\tau}, \eta_{\tau}>)$ it suffices to
consider unital completely positive maps from matrices into
$\pi_{\tau} (M_{2^{\infty}} \rtimes_{\gamma} {\Bbb Z})$.  But all of
these maps lift to unital completely positive maps into $B \subset A$.
Finally, it is clear that any abelian model defined on $\pi_{\tau}
(M_{2^{\infty}} \rtimes_{\gamma} {\Bbb Z})^{\prime\prime}$ can be
lifted to an abelian model on $A$ via the map $\Phi : A \to \pi_{\tau}
(M_{2^{\infty}} \rtimes_{\gamma} {\Bbb Z})^{\prime\prime}$.  The
desired inequality now follows easily from the definition of
CNT-entropy.  $\Box$

\begin{corollary} 
If $A$ is a unital non-type I $C^*$-algebra,  then $\infty \in
TE_{Inn} (A)$, where $TE_{Inn} (\cdot)$ denotes the inner topological
entropy invariant [BDS, 6.4.1].
\end{corollary}

{\noindent\bf Proof.}  Considering crossed products by automorphisms
with infinite CNT-entropy (for example $(\otimes_{\Bbb N}
M_{2^{\infty}}, \otimes_{n \in {\Bbb N}} \gamma^n)$) it is clear that
$\infty \in CNT_{Inn} (A)$.  The corollary then follows from [Dy,
Prop.\ 9] in the case that $A$ is exact or the remark that even the
identity automorphism has infinite topological entropy when $A$ is not
exact.  $\Box$

We conjecture that type I $C^*$-algebras are also completely
determined by the topological entropy invariant $TE_{Inn} (\cdot)$
(i.e.\ we conjecture that the topological analogue of [NS, Cor.\ 7]
is also true).  There is a natural strategy for proving this by using
a composition series to reduce to the case of continuous trace
$C^*$-algebras.  However, this would require understanding the
behavior of topological entropy in extensions and this appears to be
a highly nontrivial problem.  Another approach would be to prove that
in type I $C^*$-algebras one always has a CNT-variational principle
and hence reduce the problem to [NS, Cor.\ 7].  However, proving a
CNT-variational principle for all type I $C^*$-algebras does not
appear very easy either.

\end{document}